\documentclass[preprint,12pt]{elsarticle}

\usepackage{siunitx}
 \usepackage{graphics}
 \usepackage{subfigure} 
\usepackage[version=4]{mhchem} 
\usepackage{xcolor} 
\usepackage{float}
\usepackage{amsmath}
\usepackage{bm}
\usepackage{amssymb}
\usepackage{lineno}
\usepackage{epstopdf}
\AppendGraphicsExtensions{.tif}
\usepackage{cleveref}
\begin{document}

\begin{frontmatter}

\title{Model reduction, machine learning based global optimisation for large-scale steady state nonlinear systems}
\author[a]{Min Tao}
\author[a]{Panagiotis Petsagkourakis}
\author[a]{Jie Li}
\author[a]{Constantinos Theodoropoulos \corref{mycorrespondingauthor}}
\cortext[mycorrespondingauthor]{Corresponding author}
\address[a]{Department of Chemical Engineering , The University of Manchester, M13 9PL, UK}

\begin{abstract}
Many engineering processes can be accurately modelled using partial differential equations (PDEs), but high dimensionality and non-convexity of the resulting systems pose limitations on their efficient optimisation. In this work, a model reduction, machine-learning  methodology combining principal component analysis (PCA) and artificial neural networks (ANNs) is employed  to construct a reduced surrogate model, which can then be utilised by advanced deterministic global optimisation algorithms to compute global optimal solutions with theoretical guarantees. However, such optimisation would still be time-consuming due to the high non-convexity of the activation functions inside the reduced ANN structures. To develop a computationally-efficient optimisation framework, we propose two alternative strategies: The first one is a piecewise-affine reformulation of the nonlinear ANN activation functions, while the second one is based on deep rectifier neural networks with ReLU activation function. The performance of the proposed framework is demonstrated through two illustrative case studies.
\end{abstract}
\begin {keyword}{Model reduction, Distributed parameter systems, Piecewise affine reformulation, Deep rectifier neural networks, Data-driven optimisation, Machine learning, Global optimisation}
\end {keyword}

\end{frontmatter}

\section{Introduction}
\label{S:1}
\par Partial differential equation (PDE)-based process models, also termed distributed-parameter systems, have wide applicability in industrial engineering areas \cite{boukouvala2017global},  such as chemical \cite{tao2016hybrid}, biochemical \cite{park2015integration}, and mechanical engineering \cite{yang2008simulation} and aerodynamics \cite{kleber2001simulation}. However, complex PDEs are inherently high-dimensional and non-convex, including multiple local optima, hence resulting in intensive computational costs when  the computation of global optima is sought. Moreover, most of the generic commercial PDE simulators \cite{multiphysics1998introduction, fluent2015ansys} are essentially black-box and offer no optimisation options. Even if complex model codes are accessible in open-source software (e.g. \cite{jasak2007openfoam}), the cost of direct optimisation is often unacceptable. To date, performing optimisation tasks efficiently for large-scale complex systems, is still a challenge in engineering design. 
\par A promising way  to deal with high dimensionality  is to use  projective model order reduction  methods, which reduce the complexity of detailed models but preserve their main input-output features \cite{schilders2008model}. The popular principal component analysis (PCA) strategy, an efficient dimensionality reduction technique in  data science \cite{hinton2006reducing}, also termed as Karhunen-Loeve decomposition or proper orthogonal decomposition (POD), is usually combined with projection and/or surrogate model approaches to construct reduced models. POD together with Galerkin projection is capable of producing high-fidelity low-dimensional models for optimisation tasks \cite{theodoropoulou1998model}. Similarly, the combination of POD with machine-learing Artifical Neural Networks (ANN) has been used to construct reduced surrogate models for black-box large-scale dynamic systems, resulting in efficient optimisation and control strategies \cite{xie2015data}. In addition,  PCA and Kriging models have been utilised to  replace complex process models \cite{malik2018principal}.
\par Furthermore, \textit{equation-free} methodologies offer another effective model reduction approach for large-scale black-box systems, for optimisation and control purposes.  Exploiting the dominant eigendirections of the outputs of complex black-box system models, or direct historical system data, low-dimensional reduced Jacobian and Hessian matrixes can be computed. An equation-free based reduced SQP method was proposed exploiting the computation of low-dimensional Jacobians and Hessians, to accelerate the  optimisation procedure for large-scale steady state nonlinear systems \cite{bonis2012model}. An aggregation function was subsequently applied to address general nonlinear inequality constraints, extending the scope and capability of equation-free reduced SQP methods \cite{petsagkourakis2018reduced}.  Furthermore,  equation-free based dynamic optimisation and control methods have also been constructed \cite{bonis2013multiple,luna2005input}. An extensive discussion about model reduction based optimisation methodologies can be found in \cite{theodoropoulos2011optimisation}.
\par To address non-convexity in complex nonlinear optimisation problems, both stochastic and deterministic algorithms can be utilised. Stochastic search methods, such as simulated annealing \cite{kirkpatrick1983optimization} and genetic algorithms \cite{chambers2019practical}, can globally explore the feasible solution space avoiding local optima.  However, such stochastic search algorithms are slow for large-scale problems and offer no theoretical guarantees on the global optimality of the computed solutions. Deterministic global optimisation methods are capable of computing  global optima utilising branch-and-bound techniques \cite{floudas2005global}, but they are often computationally intensive for large-scale systems due to the need for multiple evaluations of the lower bounds of the optimisation problems. 
The aim of this work is to construct an efficient deterministic global optimisation framework for large-scale steady-state input/output (black-box) systems combining model reduction with machine-learning methodologies.  
Often a single model reduction technique cannot easily deal with the complexities of large-scale nonlinear systems. For example, although optimal principal component regressions (PCRs) \cite{pires2008selection} are popular to deal with high dimensional input-output data, the linear or low-complex models are not accurate enough to replace high nonlinear complex system models. POD on the other hand, is a very powerful method, but projecting the original system onto the global POD modes is not always easy and requires full knowledge of the full-scale system model. Meanwhile, ANN models use machine learning to capture highly nonlinear behaviours but usually require large-scale ANN structures (increasing number of neurons and layers) due to the high dimensionality of the original systems.
Combining  model reduction techniques, e.g. principal component analysis (PCA) with artificial neural networks (ANNs)\cite{lang2009reduced}, can produce accurate reduced surrogate models. Then such reduced ANN models can be explicitly utilised by global general-purpose optimisation solvers.
\par Nevertheless, performing global optimisation tasks with general ANN models is still time consuming (even for reduced ANNs), hence most existing research focuses on local optimisation and/or small-scale problems. Surrogate ANN models have been used to replace superstructure process models and have been optimised locally \cite{henao2011surrogate,fahmi2012process}. Small-scale ANN models (1 hidden layer, 3 neurons) were constructed and optimised globally by the advanced global solver BARON \cite{tawarmalani2005polyhedral}. Larger ANN models are more expensive to optimize as high non-convexity often requires the repeated use of branch-and-bound algorithms. A reduced space-based global optimisation method, recently proposed by \citet{schweidtmann2019deterministic}, projected the iteration space of non-convex variables onto the  subspace of dependent variables, resulting in small-size sub-problems and, consequently, in significant computational savings. 

In this work, two strategies are adopted to construct efficient reduced models in a PCA-ANN global optimisation framework. The first is a piecewise affine (PWA) reformulation technique of the nonlinear ANN activation function, while the second is the use of a deep rectifier neural network.
It should be noted that this work extends previous preliminary findings of the authors \cite{tao2019reduced}.

\par The rest of the paper is organized as follows. In Section 2, the basic PCA-ANN global optimisation framework is proposed and the detailed theoretical basis and implementation are provided. In Section 3, the PWA-based reformulation is outlined and illustrated with an example. In Section 4, the deep rectifier ANN-based improvement is employed in the optimisation framework and validated using a large-scale combustion case study. In Section 5, conclusions and further applications are discussed.     
\section{Problem formulation}
\label{S:2}
\par In this work, a model reduction-based optimisation framework is presented to deal with large-scale nonlinear steady-state systems focusing on the optimisation of spatially distributed processes, described by sets of dissipative PDEs of the form:
\begin{equation}
\frac{\partial \boldsymbol{y}}{\partial t}= D\{\frac{\partial  \boldsymbol{y}}{\partial x},\frac{\partial^{2} \boldsymbol {y}}{\partial x^{2}},...,\frac{\partial^{n} \boldsymbol {y}}{\partial x^{n}},\boldsymbol{d}\}+R( \boldsymbol{d},\boldsymbol{y}) \\
\end{equation}
Here $t \in \mathbb{R}$ denotes time, $x \in \mathbb{R}^{N_x}$, $N_x$ the spatial dimensions, $N_x$=1,2, or 3.  $D\in \mathbb{R}$ is the dissipative spatial differential operator, $\boldsymbol{d}\in \mathbb{R}^{N_d}$ the parameter variables  and $\boldsymbol{y}\in \mathbb{R}^{N_y}$  a set of state variables, $R(\boldsymbol{d},\boldsymbol{y}): \mathbb{R}^{N_d} \times \mathbb{R}^{N_y} \to \mathbb{R}^{N_y}$ are the nonlinear terms. Considering steady state systems and assuming that $\boldsymbol{y}(t,x)\xrightarrow{}\boldsymbol{y}(x)$, and  ${\partial \boldsymbol{y}}/{\partial t}=0$, the above equations become:
\begin{equation}
0= D\{\frac{\partial  \boldsymbol{y}}{\partial x},\frac{\partial^{2} \boldsymbol {y}}{\partial x^{2}},...,\frac{\partial^{n} \boldsymbol {y}}{\partial x^{n}},\boldsymbol{d}\}+R( \boldsymbol{d},\boldsymbol{y}) 
\end{equation}
\par  Therefore, the general optimisation problems for steady state PDE-based systems can be formulated as the following problem   $\boldsymbol{P1}$ :
  \begin{equation}
  \begin{aligned}
(\boldsymbol{P1}) \quad &\min_{\boldsymbol{d}} \quad  G(\boldsymbol{d},\boldsymbol {y})  \\   
s.t.  0&=D\{\frac{\partial  \boldsymbol{y}}{\partial x},\frac{\partial^{2} \boldsymbol {y}}{\partial x^{2}},...,\frac{\partial^{n} \boldsymbol {y}}{\partial x^{n}},\boldsymbol{d}\}+R( \boldsymbol{d},\boldsymbol{y}) \\
&\left. A \{\frac{\partial  \boldsymbol{y}}{\partial x},\frac{\partial^{2} \boldsymbol {y}}{\partial x^{2}},...,\frac{\partial^{n} \boldsymbol {y}}{\partial x^{n}}\}\right|_{x=\Omega}=h_{bds}( \boldsymbol{d},\boldsymbol{y}) \\
& g_{cons}(\boldsymbol{d},\boldsymbol {y}) \leq 0
\end{aligned}
\end{equation}
\par where $G(\boldsymbol {d},\boldsymbol{y}) : \mathbb{R}^{N_d} \times \mathbb{R}^{N_y}\to \mathbb{R}$ is the objective function. The equality constraints are the system PDEs, which in the case of a black-box system are not explicitly available, with corresponding boundary conditions.  $h_{bds}( \boldsymbol{d},\boldsymbol{y}): \mathbb{R}^{N_d} \times \mathbb{R}^{N_{y}} \to \mathbb{R}^{N_y} $ are the right hand sides of the boundary conditions, $A$ is the operator of the boundary condition equations, $\Omega$ are the boundaries and $g_{cons}( \boldsymbol{d},\boldsymbol{y}): \mathbb{R}^{N_d} \times \mathbb{R}^{N_{y}} \to \mathbb{R}^{N_y}$ denote other general constraints, e.g. bounds and other constraints, for the state variables $\boldsymbol{y}$ and the design parameter variables $\boldsymbol{d}$. 
\par In general, the unavailability of system equations inside commercial software prohibits the use of direct model-based optimisation techniques. Even in the case that large-scale system equations are available, the optimisation problem  $\boldsymbol{P1}$ can not be efficiently handled by global optimisation algorithms \cite{houska2019global}.  In this work, this barrier is overcome by employing accurate surrogate models to formulate a highly accurate approximate problem $\boldsymbol{P2}$, which is then utilised by a general purpose global optimisation solver. 

\par If we use explicit surrogate inputs-outputs to replace the black-box system equations in the above formulation, then problem $\boldsymbol{P1}$ can be transformed into the following problem $\boldsymbol{P2}$:
\begin{equation}\label{eq4}
  \begin{aligned}
(\boldsymbol{P2}) \quad \min_{\boldsymbol{d}} \quad & G(\boldsymbol {d},\boldsymbol{y'})  \\   
s.t.   \boldsymbol  {y'}&=F(\boldsymbol{d}) \\
 g_{cons}&(\boldsymbol {d},\boldsymbol{y'}) \leq 0
\end{aligned}
\end{equation}
\par where $\boldsymbol {y'}$ are the outputs of the surrogate model and $F$ the black-box nonlinear operator. 

The gap between problems $\boldsymbol{P1}$  and   $\boldsymbol{P2}$  can be measured by the output errors ($\boldsymbol {y}-\boldsymbol {y'}$). If these errors are small enough, the global optima of the reduced surrogate model-based problems $\boldsymbol{P2}$ will be close to the global solutions of the original problem $\boldsymbol{P1}$. Thus, accurate surrogate models are key  to guarantee small errors across the design domain. The accuracy of surrogate models depends on sampling quality and quantity, machine-learning techniques employed  and model building techniques. This work assumes that enough representative/informative samples are available for building accurate models. Details about efficient sampling approaches are discussed in Section 2.1.

Then model building procedures directly affect the accuracy and the complexity of constructed surrogate models, which in turn have a significant effect on the computational accuracy and speed of the subsequent deterministic global optimisation. 

This work employs a double model reduction process in conjunction with machine learning, through a combination of PCA, ANN and reformulation techniques to generate an accurate reduced model, which is subsequently used to compute near global solutions for the original problem $\boldsymbol{P1}$. In the following sections, we are discussing the basic components of our PCA-ANN-global optimisation methodology.

\subsection{Sampling and data collection }
\par  To  build accurate surrogate models, suitable sampling methods are needed to collect highly representative samples for a range of design variables. Inefficient  sampling  strategies, including too few samples and/or  unrepresentative sampling, would result in inaccurate reduced models, in turn producing inaccurate optimal solutions. While provably representative sampling is still an open problem, there are several popular sampling techniques such as Hammersley sequences \cite{faure1986star}, D‐optimal designs \cite{de1995d} and Latin Hypercube (LHC)\cite{loh1996latin} that can produce good quality results. Hammersley sequences employ quasi-deterministic sequences to convergence to a set of {\it informative} samples. 

The D-optimal design approach aims to reduce the number of experimental runs and maximise sample variances. The LHC method can produce samples covering the whole design space and maximize the difference among the generated samples. Specifically, the sample domain is divided into many sub-domains where sample points are generated randomly in order to represent the specific sub-domain. In this work, we choose LHC  because it has been shown to be able to fill the design space and to capture input/output relationships given an {\it adequate} number of samples. The number of LHC samples is decided by testing the model accuracy. In general, more LHC samples are more likely to contain the information needed to capture complex input/output relationships. LHC sampling for complex systems often requires a relatively large number of samples, which is also the pre-condition to perform successful PCA reduction and ANN-based surrogate model construction.

In the presence of constraints (such as  $g_{cons}$ here), it is hard for the LHC algorithm to directly capture the complex design space.  Previous work \cite{boukouvala2017argonaut} employed constraints to filter the LHC samples in order to reduce function evaluations for expensive systems. 
In addition, a complex strategy was used to first decompose the design space into many subdomains, where system features were represented through multiple low-fidelity models. Nevertheless the adaptive optimisation performed within the constrained sampling strategy can be computationally intensive for high-dimensional outputs $\boldsymbol{y}$ and/or inequality constraints for $\boldsymbol{y}$. 

In this work, we utilised universal ANN surrogate models to capture the nonlinear behaviour of the black-box PDE equality in its entirety. The constrained sampling strategy may possibly lead to a discontinuous design space, requiring much larger ANN structures to capture it  \cite{henao2012superstructure}. Hence, we separated the expensive black-box PDE-based equality constraints from the "known" inequality constraints.  Our aim is to construct accurate but simple ANN models to replace the PDE-based equality constraints, which together with the known inequality constraints $g_{cons}$, provide a highly accurate explicit model formulation $\boldsymbol{P2}$ to the general-purpose global optimisation solvers.

\par While building accurate surrogate models requires enough representative/informative samples, too many samples would lead to intensive computations for  high-dimensional systems. Improving sampling efficiency can significantly reduce computational times. Exploiting process knowledge or advanced adaptive sampling approaches may help to achieve this goal. Prior knowledge about the processes can provide useful information to collect representative samples with higher probability. However, this requires case-by-case detailed experience about the black-box systems. Adaptive sampling uses a few prior input/output samples to subsequently generate representative samples through solving a set of optimisation problems. Most relevant previous studies in literature \cite{eason2014adaptive,liu2018survey} deal with low-dimensional inputs/outputs, leading to small-size optimisation problems.  This work, however, deals with high-dimensional outputs $\boldsymbol{y}$, possibly together with large numbers of inequality constraints $g_{cons}$, hence the computation costs for performing adaptive sampling could be high. Nevertheless adaptive sampling procedures are fully compatible with the algorithms developed here provided that the relevant computations can be appropriately reduced.

\par Here we consider a more general approach without exploiting process knowledge and adaptive sampling techniques. The assumption is that the sampling process takes place offline and does not directly affect the computational efficiency of the online optimisation computations. Nevertheless, the proposed model reduction based global optimisation framework can be easily combined with prior knowledge and adaptive sampling approaches, as mentioned above to speed-up the offline parts of the computations. 

\par We collect samples across the space of design parameters $\boldsymbol{d}$  and corresponding input-output data sets ($\boldsymbol{D} \in \mathbb{R}^{N_d \times N}$, $\boldsymbol{Y} \in \mathbb{R}^{m\times N}$), where $m\in \mathbb{N}$ is the number of discrete interval points, which for distributed parameter systems tends to be be a large number, and $N\in \mathbb{N}$ is the number of samples.   The obtained data sets ($\boldsymbol{D}$, $\boldsymbol{Y}$)  are then used to construct accurate reduced surrogate models through the combination of PCA and ANN. 
\subsection{Principal Component  Analysis (PCA)}
\par Due to the high dimensionality of spatially discrete output data $\boldsymbol{Y}$, directly constructing surrogate ANN models would result in large ANN structures. Here, the popular PCA method is first employed to build a reduced model from output data $\boldsymbol{Y}$.  
\par A sampling method (here LHC as discussed in the previous section) is firstly employed to construct a data ensemble $\boldsymbol Y$over a finite spatial interval $\Omega' \in \mathbb{R}$.  PCA  then calculates a reduced output $\boldsymbol {U}$,  by projecting the data sample $\bf Y$ onto the subspace of the "small" set of principal components (PCs)  $P= (p_1,p_2,...,p_k)$, $k\in \mathbb{N}$ being the number of PCs.
\begin{equation}\label{eq5}
\begin{aligned}
\boldsymbol{U}=\boldsymbol{PY}
  \end{aligned}
\end{equation}
Here $\boldsymbol{U}\in \mathbb{R}^{k\times N}$ is the projection of the original data $\boldsymbol{Y}$ onto the reduced subspace $\mathbb {P}$ and $\boldsymbol{P}\in \mathbb{R}^{k\times m}$ is the corresponding orthogonal projector. 
In the PCA method the matric $\boldsymbol {P}$ is constructed through the covariance matrix, $\boldsymbol{C_y}\in \mathbb{R}^{m\times m}$ of the output data $\boldsymbol{Y}$:
\begin{equation}
\begin{aligned}
\boldsymbol{C_y}=\frac{1}{m-1}\boldsymbol{YY^T}
  \end{aligned}
\end{equation}
Here we seek to minimise covariance between data and maximise variance i.e. minimise the off-diagonal elements of $\boldsymbol{C_y}$, while maximising its diagonal elements. This is equivalent to performing singular value decomposition (SVD) on $\boldsymbol{C_y}$:  

\begin{equation}\label{eq7}
\begin{aligned}
\boldsymbol{C_y}=\boldsymbol{Z^TZ}=(\frac{1}{\sqrt{m-1}}\boldsymbol{Y}^T)^T(\frac{1}{\sqrt{m-1}}\boldsymbol{Y}^T)=\boldsymbol{VDV^T}
  \end{aligned}
\end{equation}
where $D \in \mathbb{R}^{m\times m}$ is a diagonal matrix whose diagonal elements are the eigenvalues of $Z^TZ$ and $V$ is the orthogonal matrix whose columns are the eigenvectors of $Z^TZ$, which as can be easily shown are equivalent to the principal components of $\boldsymbol Y$. In fact we can keep the first $k$ PCs corresponding to the $k$ dominant eigenvalues of $\boldsymbol {C_y}$, where usually $k<<m$, hence $\boldsymbol V \in \mathbb{R}^{m\times k}$ and $D$ now contains only the $k$ most dominant eigenvalues of the system, $D\in \mathbb{R}^{k\times k}$. We can then set $\boldsymbol {P}=\boldsymbol{V}^T$ and perform data reduction through the projection in Eq. \ref{eq5}.
The original data sample, $\boldsymbol Y$ can be reconstructed from the projected data:
\begin{equation}\label{eq71}
\begin{aligned}
\boldsymbol{Y}=\boldsymbol{P^T}\boldsymbol{U}
  \end{aligned}
\end{equation}More details about the theory and application of PCA can be found in \cite{richardson2009principal,hotelling1933analysis,berkooz1993proper,jirsa1994theoretical,park2008electro}.

\par  The PCA step aims to project the high dimensional output states arising from the discretisation of the PDE system, onto a small number of dominant variables. Then the resulting low-dimensional relationship between inputs and projected outputs can be captured through small-size ANN structures, hence producing efficient doubly reduced models that will significantly reduce deterministic global optimisation computations. Implementing the PCA step is not always easy, as it is very sensitive to the quantity and quality of samples. Only with enough representative/informative samples can PCA be efficiently performed and globally capture the accurate dominant PCs. 

\subsection{Artificial Neural Networks (ANNs)}
\par We employ ANNs on the reduced models (Eq.\ref{eq5}) from the PCA step. ANN-based models are chosen due to both successful practices and proven theoretical support that a {\it shallow} feed-forward neural network with one single layer is sufficient to represent any smooth function \cite{hornik1989multilayer}. Furthermore, advanced optimisation algorithms have been developed to handle the machine-learning component of the process, i.e. handling the the manipulated variables for ANN structures, such as Levenberg-Marquardt backpropagation \cite{yu2011levenberg} and Bayesian regularization backpropagation \cite{mackay1992practical}. Fig.\ref{fig1} shows a conventional feed forward neural network with a hyperbolic 
tangent activation function $tanh(\cdot)$.
\textbf{\begin{figure}[H] 
    \centering
    \subfigure{
\includegraphics[width=0.4\textwidth]{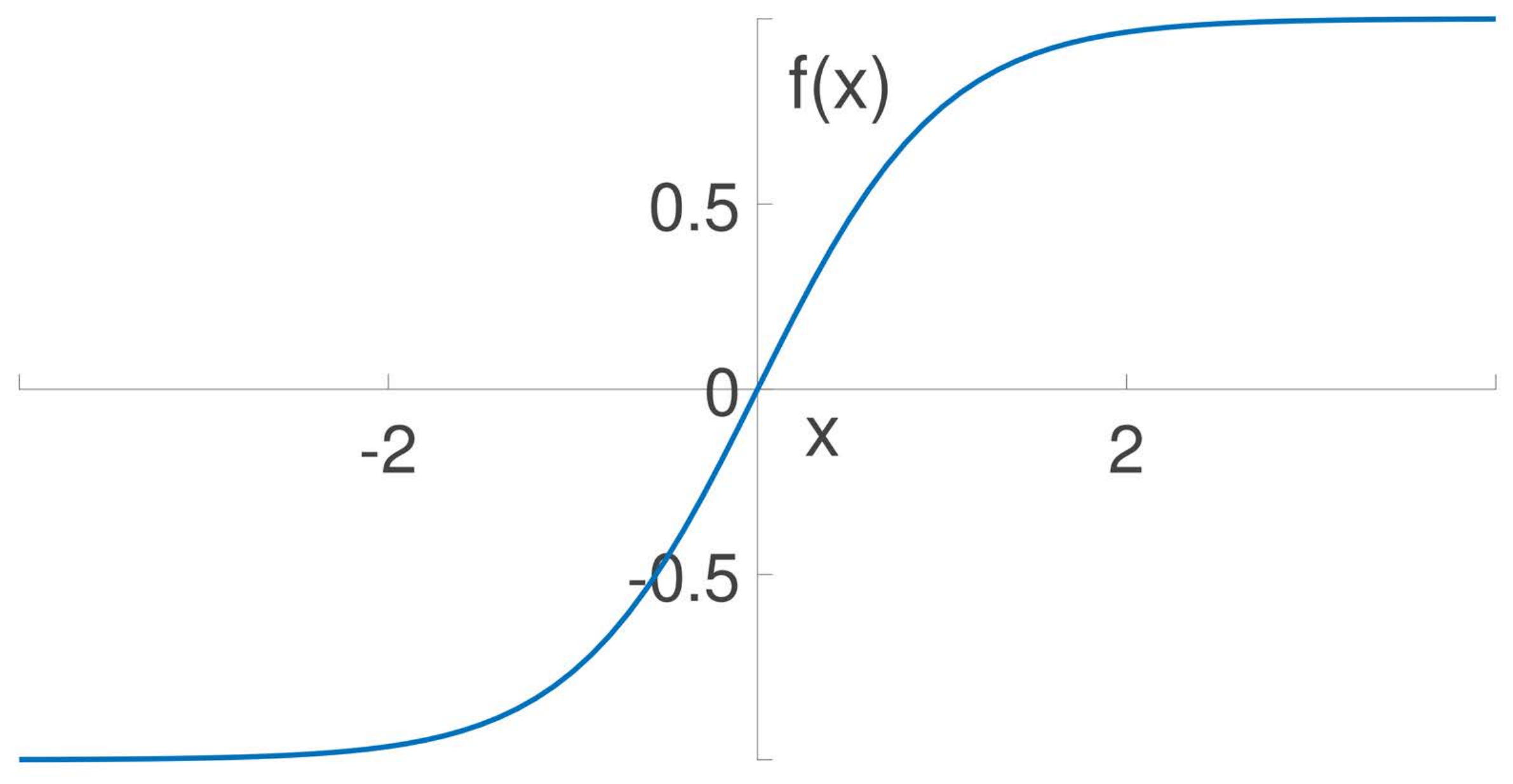}}
\hfill
\subfigure{
    \includegraphics[width=0.5\textwidth]{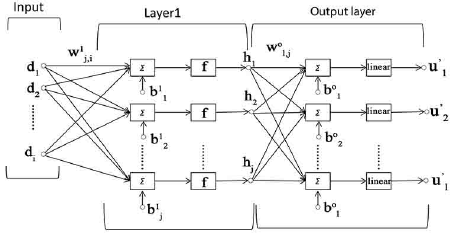}}
    \caption{Feed-forward neural network with hyperbolic tangent activation function}
     \label{fig1}
\end{figure}}
Shallow ANNs, as the one displayed in Fig.\ref{fig1}, are implemented in our basic PCA-ANN global optimisation framework. The feed-forward ANN contains three main components: The input layer, the hidden layer (only one in a shallow ANN) and the output layer, which sequentially perform transformations on the input variables. The input variables, $\boldsymbol{d}=( d_1,d_2,...,d_{N_d})$, are first linearly transformed and then non-linearly activated through the hidden layer, and further forced by linear transformation and sequential activation in the output layer, to finally formulate the output variables $\boldsymbol{u'}=( u'_1,u'_2,...,u'_{k}), \in \mathbb{R}^{k}$. The mathematical description is given in Eq.\ref{eq17} : 
\begin{equation}\label{eq17}
\begin{aligned}
h_j&=f(\sum_{i=1}^{N_d}w^1_{j,i}d_i+b^1_j), &\forall j \in \lbrace 1,2,...,n \rbrace\\
u'_l&= \sum_{j=1}^nw^o_{l,j}h_j+b^o_l, &\forall l \in \lbrace 1,2,...,k \rbrace
  \end{aligned}
\end{equation}
Here $h_j\in \mathbb{R}$ is the output value from the hidden layer with $n\in \mathbb{N}$ neurons,  $j=1,\dots, n$ and $f\in \mathbb{R}$ is the activation function. Each neuron $j$ contains two parameters: weights $w^1_{j,i}\in \mathbb{R}$ and biases $b^1_j\in \mathbb{R}$ which perform linear transformations. Similarly, $u'_l\in \mathbb{R}$ is the final value from the output layer with $k$ neurons, $l=1,\dots,k$, including weights $w^o_{l,j}\in \mathbb{R}$ and biases $b^o_l\in \mathbb{R}$. Three activation functions, the sigmoidal, the hyperbolic tangent and the linear function, are widely used in neural networks. In this work, the hyperbolic tangent function $f$ was utilised to convert the output value into the range [0,1] in the hidden layer while the linear function was applied in the output layer. 
The configured feed-forward neural network was subsequently trained through the back-propagation algorithm using the reduced low-dimensional data sets ($\boldsymbol{D}$, $\boldsymbol{U}$) from the PCA step.  
Similar to the PCA technique, the performance of surrogate ANN models is significantly affected by quantity and quality of collected samples. 

\par The above ANN training can be viewed as an application of supervised machine learning methodology \cite{hastie2009elements}, which aims to learn a function that maps the  input variables $\boldsymbol{d}$ to the output variables $\boldsymbol{u'}$  based on the collected input-output samples  ($\boldsymbol{D}$, $\boldsymbol{U}$). Here the structure of the learned function is the chosen ANN model with the corresponding learning algorithm being the training method, the  Levenberg-Marquardt algorithm with an early stopping procedure. Performing this machine learning task will generate a predictive ANN model, which together with the previous PCA projection will formulate a low-dimensional surrogate model for the original high-dimensional system.

\subsection{PCA-ANN global optimisation framework}
\par
\begin{figure}[H]
    \centering
    \includegraphics[width=0.4\textwidth]{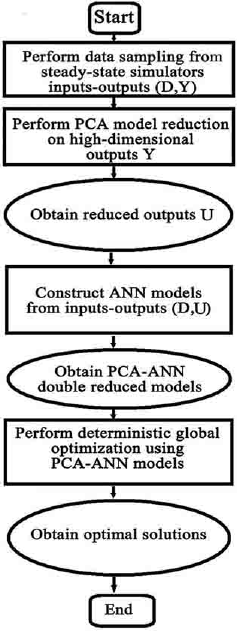}
    \caption{Flow chart of the basic PCA-ANN global optimisation framework}
    \label{fig2}
\end{figure}

 To cope with the non-convexity of highly non-linear systems, deterministic optimisation methods are considered for the reduced surrogate model from the PCA-ANN reduction. The black-box or grey-box  global optimisation problem can be transformed into the general explicit NLP optimisation problem as follows combining Eqs.(\ref{eq4},\ref{eq5},\ref{eq17}):
\begin{equation}\label{eq18}
\begin{aligned}
\min_{\boldsymbol{d}=[d_1,d_2,...,d_{N_{d}}]} \quad & G(\boldsymbol{d},\boldsymbol{u'})   \\   
s.t.   h_j&=f(\sum_{i=1}^{N_d}w^1_{j,i}d_i+b^1_j), &\forall j \in \lbrace 1,2,...,n \rbrace\\
u'_{l}&= \sum_{j=1}^nw^o_{l,j}h_j+b^o_l, &\forall l \in \lbrace 1,2,...,k \rbrace \\
\boldsymbol{u'}&=( u'_1,u'_2,...,u'_{k}),\\
 g_{cons}&(\boldsymbol{d},\boldsymbol{P^T}\boldsymbol{u'}) \leq 0 
\end{aligned}
\end{equation}
In this work, $g_{cons}(\boldsymbol{d},\boldsymbol{P^T}\boldsymbol{u'})$ includes the box bound constraints and possible inequality constraints for  design variables $\boldsymbol{d}$ and discretised state variables $\boldsymbol {y'}$.  The ANN-based nonlinear objective function $G$, can be reformulated into the constraints.  The main non-convexity of the optimisation problems lies on the surrogate model constraints $h_j=f(\cdot)$ due to the highly non-convex activation function, here the hyperbolic tangent function $tanh(\cdot)$ in the feed-forward ANN structure. General-purpose global optimisation commercial software, including ANTIGONE \cite{misener2014antigone}, BARON \cite{tawarmalani2005polyhedral} and SCIP \cite{rehfeldt2018scip}, are efficient tools for the above problems due to the advanced bound tightening and branching techniques. Nevertheless, these general global solvers can not handle the $tanh(\cdot)$ formulation directly, as high performance algorithms need the explicit model equations. Therefore the explicit algebraic form $tanh(z)={(e^z-1)}/{(e^z+1)}$ is required \cite{smith2012cfd}. The basic formulation is further transformed into $tanh(z)= {-2}/{(e^z+1)}+1$ in order to produce a tighter under-estimator for the global solver \cite{schweidtmann2019deterministic}. The flow chart of the basic PCA-ANN global optimisation framework is shown in Fig.\ref{fig2}. 

\section{Piece-wise affine based formulation}
In this section,  a piecewise affine (PWA) reformulation is introduced to deal with the non-convex hyperbolic tangent activation function in the reduced surrogate ANN-based model. Previous research has suggested the PWA technique for ANN models \cite{henao2012superstructure}, which has been verified to be efficient \cite{kessler2017efficient}. 
Although these studies provided some computational results, further detailed implementation schemes and analysis have not been reported. In this work, the PWA reformulation was utilised to approximate the highly non-convex NLP problem with a MILP problem. The global optimisation algorithms for both  NLP and  MILP problems are based on the branch and bound framework. However, the branching step is performed on continuous variables for the NLP problems and on  auxiliary binary variables for the MILP problems through the use of CPLEX 12.7.1. 
An adaptive procedure to construct PWA models is presented below. 

\subsection{Adaptive procedure}
The hyperbolic tangent activation function $f(z)=tanh(z)$ is an odd function with central symmetry, which is concave on $(0,+\infty]$ and convex on $[-\infty,0)$. Therefore the PWA approximation on $[-\infty,0)$ can be directly computed from the PWA formulation on $(0,+\infty]$. Within the range of $(0,+\infty]$, $tanh(z)$ function first increases and then tends to level off with a slight increase towards the limit value of 1. 
The adaptive PWA procedure starts from the interval $(0,+\infty]$ and two points, the
point of symmetry and one point close to the maximum value (equal to 1). 
Then a new point is chosen between the two original points so that the error $E_{ero}$ between 
$f(z)$ and its PWA approximation $f_{PWA}(z)$ (currently consisting of two intervals) is minimised.
\begin{equation}\label{eq20}
\begin{aligned}
E_{ero}&=\int abs (f(z)-f_{PWA}(z)) {\rm  d}z, 
  \end{aligned}
\end{equation}
Then the  segment with the largest error is chosen and a new point is added within to minimise $E_{ero}$ in this segment. This procedure continues iteratively until the error in eq. \ref{eq20}
becomes less than a pre-defined tolerance. Finally the points chosen for the $(0,+\infty]$ interval are mirrored to the $[-\infty,0)$ interval. 

\begin{figure}[H]
    \centering
    \includegraphics[width=0.8\textwidth]{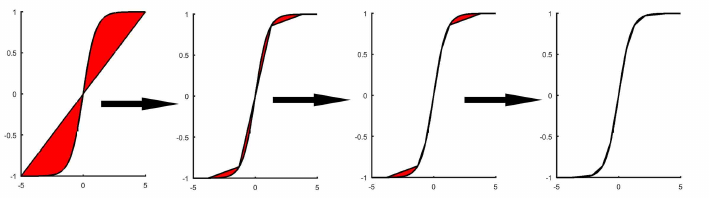}
    \caption{Adaptive PWA procedure for the hyperbolic activation function}
    \label{adap_pro}
\end{figure}

The iteration procedure efficiently produces a tight PWA representation of the $tanh(z)$ function. Fig.\ref{adap_pro} shows the adaptive process, narrowing the interval sizes and reducing the error (red shade) between piecewise affine models and $tanh(z)$. 
There are different approaches to formulate the PWA models, such as the classic method, the linear segmentation method,the  convex hull method and special structure methods. The classic method is the basic step of the other three methods and may work less efficiently \cite{misener2009global}. Here we employed the special order sets technique since the advanced MILP solvers including CPLEX could smartly exploit the structures of special order sets and speed up computations \cite{misener2010piecewise}. For $N'+1$ generated grid points $z_1,z_2,...,z_{N'+1} \in \mathbb R$ and correspondingly $N'$ linear models, the general PWA formulation introducing the special sets variables $h'_i$ and $\lambda'_i$, is as follows \cite{floudas1995nonlinear}:
\begin{equation}\label{eq21}
\begin{aligned}
f(z)&\approx f_{PWA}(z)= \sum_{i=1}^{N'+1}\lambda'_i f(z_i),\\
z&= \sum_{i=1}^{N'+1}\lambda'_i z_i,\\
\sum_{i=1}^{N'+1}&\lambda'_i =1,\\
\lambda'_1&\leq h'_1,\\
\lambda'_i&\leq h'_i +h'_{i-1} , \forall i \in \lbrace 2,3,...,N'\rbrace\\
\lambda'_{N'+1}&\leq h'_{N'},\\
\lambda'_i&\geq0, \forall i \in \lbrace 1,2,...,N'+1 \rbrace\\
\sum_{i=1}^{N'}&h'_i =1,\\
h'_i &\in \lbrace{0,1}\rbrace^{N'}
\end{aligned}
\end{equation}
It should be noted here that the above formulation allows only two adjacent $\lambda_i$'s to be non-zero.

Substituting the highly non-convex $f(z)$ in the PCA-ANN optimisation formulation (Eq.\ref{eq18})with the above PWA reformulation (Eq.\ref{eq21}), the general PCA-ANN-PWA based MILP optimisation problem can be obtained:
\begin{equation}\label{eq22}
\begin{aligned}
\min_{\boldsymbol{d}} \quad & G(\boldsymbol {d},\boldsymbol{u'})   \\   
s.t.   z^j&=\sum_{i=1}^{N_d}w^1_{j,i}d_i+b^1_j, &\forall j \in \lbrace 1,2,...,n \rbrace\\
h_j&= \sum_{i=1}^{N'+1}\lambda^j_i f(z_i), &\forall j \in \lbrace 1,2,...,n \rbrace\\
z^j&= \sum_{i=1}^{N'+1}\lambda^j_i z_i, &\forall j \in \lbrace 1,2,...,n \rbrace\\
\sum_{i=1}^{N'+1}&\lambda^j_i =1,&\forall j \in \lbrace 1,2,...,n \rbrace\\
\lambda^j_1&\leq h^j_1,&\forall j \in \lbrace 1,2,...,n \rbrace\\
\lambda^j_i&\leq h^j_i +h^j_{i-1} , &\forall i \in \lbrace 2,3,...,N' \rbrace, \forall j \in \lbrace 1,2,...,n \rbrace\\
\lambda^j_{N'+1}&\leq h^j_{N'},&\forall j \in \lbrace 1,2,...,n \rbrace\\
\lambda^j_i&\geq0, &\forall i \in \lbrace 1,2,...,N'+1 \rbrace, \forall j \in \lbrace 1,2,...,n \rbrace\\
\sum_{i=1}^{N'}&h^j_i =1&\forall j \in \lbrace 1,2,...,n \rbrace\\
h^j_i &\in \lbrace{0,1}\rbrace^{N'}&\forall j \in \lbrace 1,2,...,n \rbrace\\
u'_{l}&= \sum_{j=1}^nw^o_{l,j}h_j+b^o_l, &\forall l \in \lbrace 1,2,...,k \rbrace \\
\boldsymbol{u'}&=( u'_1,u'_2,...,u'_{k}),\\
 g'_{cons}&(\boldsymbol{d},\boldsymbol{P^T}\boldsymbol{u'}) \leq 0 \\
\end{aligned}
\end{equation}where $g'_{cons}$ is a PWA formulation of the possibly nonlinear inequality constraints. The following case studies consider optimisation problems with linear inequality constraints, which is enough to validate the efficiency of PWA formulation of nonlinear ANN models on computational cost and accuracy.

\section{Deep rectifier neural network based formulation}
\par The  PCA-ANN-PWA global optimisation framework may work efficiently for some applications. However, the PWA step will in general lead to additional approximation error especially for large-scale problems.To preserve the computational accuracy and still use the advanced MILP solver, the continuous piece-wise linear activation function is introduced and directly embedded in the ANN structures. 
Past efforts in computer science have developed efficient activation functions, such as the sigmoid and the $tanh(\cdot)$ function. The S-shaped sigmoid function can transfer any input signal into the range [0,1] while  the zero centered $tanh(\cdot)$ function can map the output values in the interval [-1,1]. Both of them can learn features of high nonlinear functions efficiently.
Nevertheless, the high non-convexity of these functions makes ANN training hard.  The continuous piecewise linear functions, including the $ReLU$ function and its variants, have been adopted to deal with this problem. In this work, the widely applied  $ReLU$ function is utilised in order to reserve computational accuracy and maintain the ability to use the advanced MILP solver. Nevertheless, shallow neural networks require a significantly larger number of nodes in one hidden layer to successfully represent a complex function, while deep neural networks result in more complex and non-convex training errors \cite{lee2018machine} due to their multi-layer  structure. Low-complexity two- or three-hidden layer NNs are, however, enough to capture the low-dimensional nonlinear behaviour of PCA-reduced systems.

\begin{figure}[H]
\centering
\subfigure[ReLU activation function]{
\includegraphics[width=0.4\textwidth]{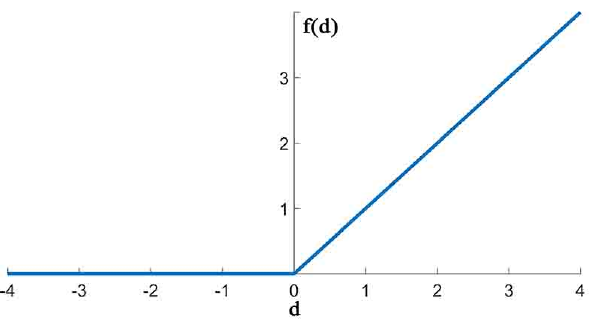}}
\hfill
\subfigure[Deep artificial neural network]{
\includegraphics[width=0.5\textwidth]{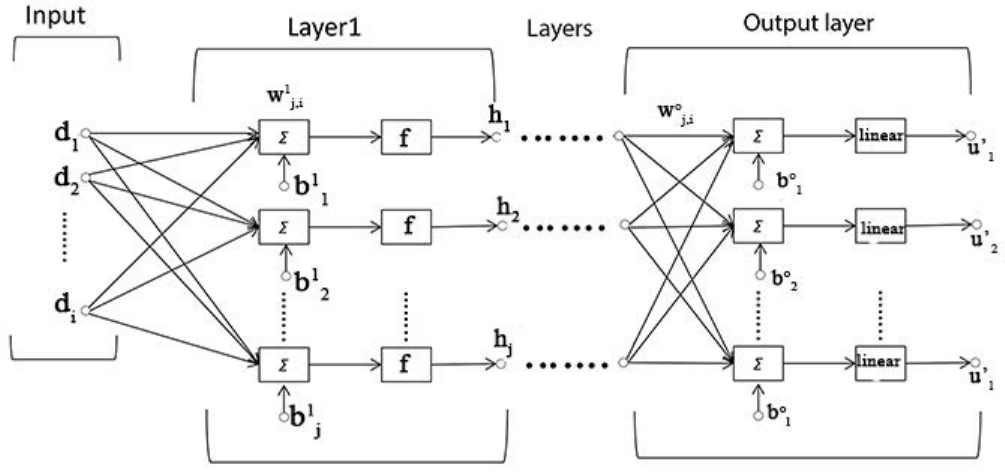}}
\caption{Deep neural network with rectified linear units}
\label{fig:8}
\end{figure}

Although deep rectifier NN-based MILP problems have been formulated in previous studies \cite{grimstad2019relu},  the combination of PCA and deep rectifier NN has not been reported before to the best of our knowledge. The mathematical equations for the deep rectifier NNs are similar to those in Fig.\ref{fig1} with more hidden layers and activation function $f(z)=max(0,z)$, which can be reformulated into a piecewise linear function through the big-M method \cite{belotti2010disjunctive}:
 \begin{equation}\label{eq24}
\begin{aligned}
z^{j_1}_1&=\sum_{i=1}^{N_d}w^1_{j_1,i}d_i+b^1_{j_1}, &\forall {j_1} \in \lbrace 1,2,...,n_1 \rbrace\\
z^{j_1}_1&={z'}^{j_1}_1-{z''}^{j_1}_1, &\forall {j_1} \in \lbrace 1,2,...,n_1 \rbrace\\
{z'}^{j_1}_1&\leq M_1(1-{bz}^{j_1}_1), &\forall {j_1} \in \lbrace 1,2,...,n_1 \rbrace\\
{z''}^{j_1}_1&\leq M_1{bz}^{j_1}_1, &\forall {j_1} \in \lbrace 1,2,...,n_1 \rbrace\\
z^{j_2}_2&=\sum_{j_1=1}^{n_1}w^2_{j_2,j_1}{z'}^{j_1}_1+b^2_{j_2}, &\forall j_2 \in \lbrace 1,2,...,n_2 \rbrace\\
z^{j_2}_2&={z'}^{j_2}_2-{z''}^{j_2}_2, &\forall {j_2} \in \lbrace 1,2,...,n_2 \rbrace\\
{z'}^{j_2}_2&\leq M_2(1-{bz}^{j_2}_2), &\forall {j_2} \in \lbrace 1,2,...,n_2 \rbrace\\
{z''}^{j_2}_2&\leq M_2{bz}^{j_2}_2, &\forall {j_2} \in \lbrace 1,2,...,n_2 \rbrace\\
&\cdot  &\cdot\\
&\cdot &\cdot \\
&\cdot &\cdot\\
z^{j_{\theta}}_{\theta}&=\sum_{j_{{\theta}-1}=1}^{n_{{\theta}-1}}w^{\theta}_{j_{\theta},j_{{\theta}-1}}h^{j_{{\theta}-1}}_{{\theta}-1}+b^{\theta}_{j_{\theta}}, &\forall j_{\theta} \in \lbrace 1,2,...,n_{\theta} \rbrace\\
z^{j_{\theta}}_{\theta}&={z'}^{j_{\theta}}_{\theta}-{z''}^{j_{\theta}}_{\theta}, &\forall {j_{\theta}} \in \lbrace 1,2,...,n_{\theta} \rbrace\\
{z'}^{j_{\theta}}_{\theta}&\leq M_{\theta}(1-{bz}^{j_{\theta}}_{\theta}), &\forall {j_{\theta}} \in \lbrace 1,2,...,n_{\theta} \rbrace\\
{z''}^{j_{\theta}}_{\theta}&\leq M_{\theta}{bz}^{j_{\theta}}_{\theta}, &\forall {j_{\theta}} \in \lbrace 1,2,...,n_{\theta} \rbrace\\
{z'}^{j_i}_i&\geq 0, &\forall {i} \in \lbrace 1,2,...,{\theta} \rbrace, \forall j_{i} \in \lbrace 1,2,...,n_i \rbrace\\
{z''}^{j_i}_i&\geq 0 &\forall {i} \in \lbrace 1,2,...,{\theta} \rbrace, \forall j_{i} \in \lbrace 1,2,...,n_i \rbrace\\
{bz}^{j_i}_i&\in \lbrace{0,1}\rbrace&\forall {i} \in \lbrace 1,2,...,{\theta} \rbrace, \forall j_{i} \in \lbrace 1,2,...,n_i \rbrace\\
u'_{l}&= \sum_{j_{\theta}=1}^{n_{\theta}}w^o_{l,j_{\theta}}{z'}^{j_{\theta}}_{\theta}+b^o_{l}, &\forall l \in \lbrace 1,2,...,k \rbrace
  \end{aligned}
\end{equation}
where $M_{i}$ is the big-M constant, 
here equal to 10000, which was enough to capture the global optimum of the subsequent tubular reactor example without numerical issues. ${z'}^{j_i}_i$ and ${z''}^{j_i}_i$ are  the auxiliary non-negative variables, ${bz}^{j_i}_i$ is the auxiliary binary variable and $h^{j_i}_i$ is the output value from the  $j_{i}$th ReLU based neuron of the $i$th hidden layer. $\theta$  is the number of hidden layers and $n_i$ is number of neurons at the $i$th hidden layer.
 Substituting the ANN model equations in the PCA-ANN optimisation formulation (Eq.\ref{eq18}) with the above reformulation (Eq.\ref{eq25}), the following PCA-DNN(ReLU) based MILP optimisation formulation can be obtained:
\begin{equation}\label{eq25}
\begin{aligned}
\min_{\boldsymbol{d}=[d_1,d_2,...,d_{N_{d}}]} \quad & G(\boldsymbol {d},\boldsymbol{u'})   \\   
s.t.  \quad  &Eq.(\ref{eq24}); \\
\boldsymbol{u'}&=( u'_1,u'_2,...,u'_{k}),\\
 g'_{cons}&(\boldsymbol {d},\boldsymbol{P^T}\boldsymbol{u'}) \leq 0 \\
  \end{aligned}
\end{equation}
where Eq.\ref{eq24} denotes the equality constraints constructed using deep neural networks through the big-M method.
This way, an improved framework is  formulated using a deep neural network (DNN) with rectified linear units (ReLU) as illustrated in Fig.\ref{fig:8}, which is first tested with the tubular reactor, and then with a more challenging large-scale combustion process. 

Moreover, it should be noticed that although ReLU activation function is the most widely applied activation function for the hidden layers of deep neural networks (DNN), sigmoid and tanh activation functions are still applied for the output layers of DNN and in RNNs, respectively \cite{nwankpa2018activation}. Different activation functions have both advantages and dis-advantages of their own depending on the type of system that we design \cite{sharma2017activation}. For example, the sigmoid function  appears in the output layers of DL architectures, which are used for predicting probability based outputs and have been applied successfully in binary classification problems, modelling logistic regression tasks as well as other neural network domains. While the tanh functions have been used mostly in recurrent neural networks for natural language processing and speech recognition tasks. The main advantage of using the rectified linear units in computation \cite{sharma2017activation}, is that, they guarantee faster computation since there is no need to compute exponentials and divisions, with overall speed of computation enhanced.  In some cases, the tanh based NNs work better than ReLU based NNs to capture the nonlinear behaviour of systems \cite{mou2017deep}. The proposed PWA approach could be employed in both sigmoid and tanh activation function. For the nonlinear systems which the Neural networks with sigmoid (maybe in output layer) and tanh activation functions can accurately represent, PWA approach still can speed up of the optimisation process for the trained NNs. 

\section{Case study 1: tubular reactor}

Here a typical tubular reactor  where an exothermic reaction takes place \cite{jensen1982bifurcation}, was used as an illustrative example to show how quantity and quality of samples affect PCA projection errors.  
The reactor model consists of 2 differential equations in dimensionless form: 

\begin{equation}\label{eq23_1}
\begin{aligned}
&0=\frac{1}{Pe_1}\frac{\partial^{2}C}{\partial y^{2}}-\frac{\partial C}{\partial y}+Da(1-C)exp(T/(1+T/\gamma)) \\
&0=\frac{1}{LePe_2}\frac{\partial^{2}T}{\partial y^{2}}-\frac{1}{Le}\frac{\partial T}{\partial y}-\frac{\beta}{Le}T+BDa(1-C)exp(T/(1+T/\gamma))+\frac{\beta}{Le}T_w   \\  
b.c. &\\
&\frac{\partial C}{\partial y}-{Pe_1}{C} =0, \frac{\partial T}{\partial y}-{Pe_2}{T} =0, \quad at \quad y=0\\
&\frac{\partial C}{\partial y} =0, \frac{\partial T}{\partial y} =0, \quad at \quad y=1
\end{aligned}
\end{equation}
$C$ and $T$ are the dimensionless concentration and temperature respectively, while $C_{exit}$ is dimensionless output concentration. $Da$ is the Damköhler number, $Le=1$ the Lewis number, $Pe_1=5$ the Peclet number for mass transport and  
$Pe_2=5$ for heat transport, $\beta=1.5$ a dimensionless heat transfer coefficient, $C$ the dimensionless adiabatic  temperature rise, $\gamma=10$ the dimensionless activation energy, $T_w=0$ the adiabatic wall temperature and $y$ the dimensionless longitudinal coordinate. The model equations were discretised using central finite differences over 250 computational nodes, resulting in 
500 algebraic equations, which comprise our in-house simulator, which was subsequently treated as a black-box (input/output) code. The initial aim was then to represent the nonlinear behaviour of 500 outputs (distributed dimensionless temperature and concentration) with respect to the single design variable Da, varying in the range [0.121, 0.400].  
For comparison purposes, we generated 6 different sample groups of
10, 20, 30, 40, 50, 60 LHC samples, respectively. Then PCA projection was employed to compute the dominant PCs, with the maximum energy/variance ratio  set to 99.8 \%. 
Finally, we utilised 500 uniform design points in the above range [0.121, 0.400], resulting in a total of 250000 outputs, to test the PCA model prediction accuracy outside the sample set utilised to generate the corresponding PCs. The results are shown in Tb.\ref{projection accuracy}.
\begin{table}[H]
\centering
  \caption{Comparative results of PCA  for different numbers of samples}
  \label{projection accuracy}
  \smallskip
  \begin{tabular}{c|c|c|c}
    \hline
 No of Samples&No of PCs&Total rel. error (250000 points)&Max error\\
  \hline
10&2&5220&3.16\\
20&2&5350&3.23\\
30&2&3680&0.86\\
40&2&4060&3.03\\
50&3&1170&0.16\\
60&3&1210&0.50\\
70&3&1077&0.49\\
80&3&1111&0.24\\
  \hline
  \end{tabular}
\end{table}
\par Tb.\ref{projection accuracy} shows that the smaller sample groups (10, 20, 30, 40, LHC samples respectively) require only 2 PCs to capture the variance of the corresponding data sets while the larger sample groups  (50, 60 LHC samples, respectively) need 3 PCs, indicating that fewer samples would miss some global information (PCs) of the whole design space. This is also the possible reason that the resulting PCA projections from the smaller sample groups generate higher errors in the  validation process.  One exception to this is the 40-sample case, which is less accurate than 30-sample one  despite the higher number of samples. The possible explanation is that the generated 40 samples happen to be less representative than the generated 30 LHC samples. 
%
\par Importantly, the PCA step reduces the original 500 distributed state variables into only 3 PC variables, which significantly decreases the computations for training ANN models and for subsequently performing deterministic global optimisation of the trained models. We employed the 50 LHC samples group for training our ANNs. To avoid over- and under-fitting, the defined domain was randomly divided into a training, a validation and a test set with respective size ratios of 0.7 : 0.15 : 0.15.  The machine-learning algorithm was implemented through the MATLAB Neural Network Toolbox to fit the corresponding weights and biases by minimizing the mean squared error (MSE) between the ANN model and the training set using Levenberg-Marquardt algorithm and the early stopping procedure. To obtain a suitable number of neurons in the hidden layer, the training process was repeated using an increasing number of neurons until the MSE for all three sets became less than a pre-defined tolerance, here $1\times 10^{-4}$. We tested the training process for a 5-neuron shallow ANN generated for the PCA-reduced model with a single input and three outputs, which took less than 1s for the 849 iterations required for convergence. Without the PCA reduction step, each iteration requires 12.41s for the single input and 500 output ANN with the same hidden layer structure,  to slowly reduce the training errors. Nevertheless, convergence could not be achieved in 20 consecutive training runs for this ANN due to large validation and/or test errors, implying that a larger ANN with more layers is needed to accommodate the relevant input/output information. This would in turn lead to a much higher number of computations for deterministic global optimisation of the trained ANN. All runs were performed  on a Desktop (Intel® Core(TM) CPU 2.5 GHz, 16 GB memory, 64-bit operating system) running Windows 10. 

\par We subsequently investigated the computational efficiency of the PCA-ANN-PWA and PCA-DNN optimisation framework through the maximisation of the exit concentration of the tubular reactor system by manipulating the temperature of 3 cooling zones along its length. The general mathematical formulation of the optimisation problem is as follows :
\begin{equation}\label{eq23}
\begin{aligned}
&\qquad \qquad  \quad  \qquad\max \limits_{T_{wi}}C_{exit}\\
s.t.& \qquad  PCA-ANN \qquad model \\
 &T_w(y)=\sum_{i=1}^{3}\left(H(y-y_{i-1})-H(y)-y_i\right)T_{wi} \\
 & 0 \leq T_{wi}\leq 4 
\end{aligned}
\end{equation}
Here 

$C_{exit}$ is the dimensionless output concentration. The system parameters are $Pe_1$ = 5, $Pe_2=5$, $Le$ = 1, $\beta$ = 1.5, $\gamma $  = 20, $B$ = 12, $D_a$=0.1.  $T_w$ is  the adiabatic wall temperature, $T_{wi}$ are the corresponding temperatures at the three cooling zones, $i$ and $H$ is the Heaviside step function.

\begin{figure}[H]
\centering
\subfigure[Solution profiles for temperature ]{
\includegraphics[width=0.4\textwidth]{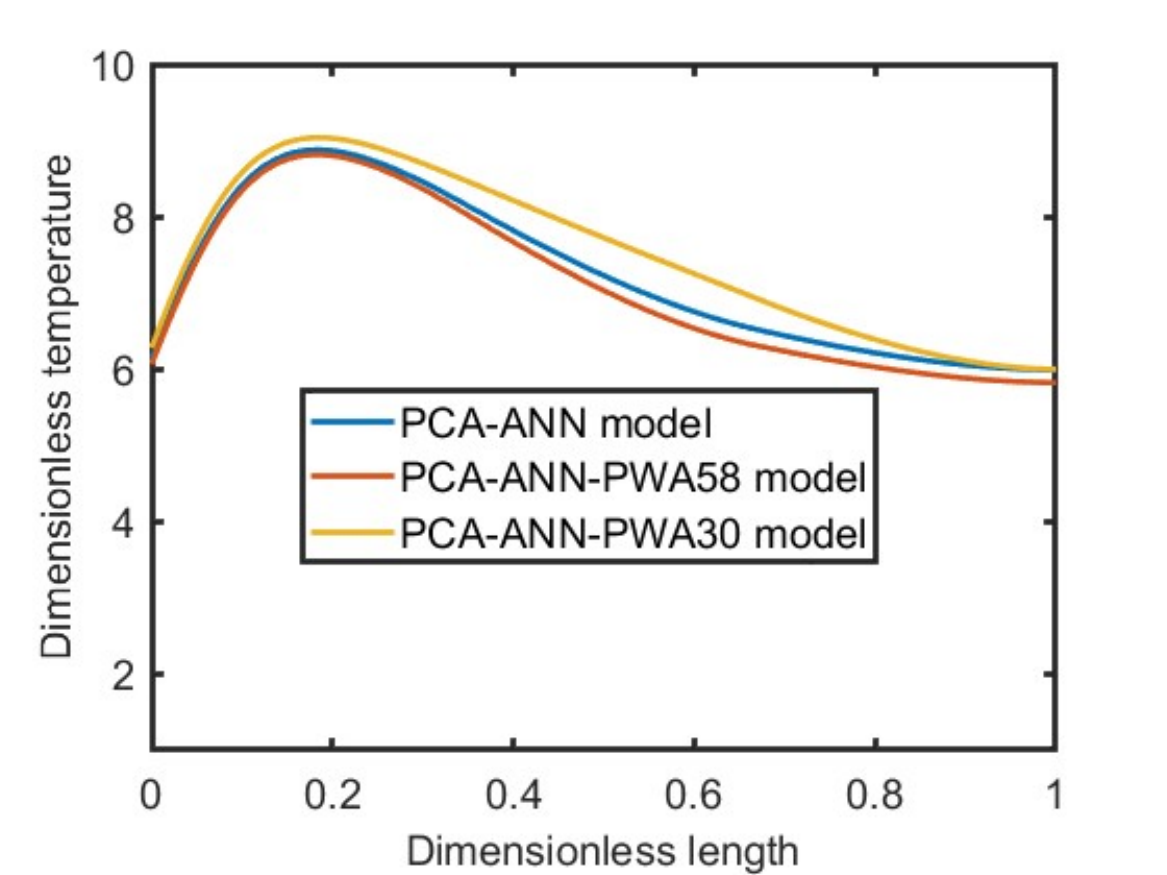}}
\hfill
\subfigure[Solution profiles for concentration]{
\includegraphics[width=0.4\textwidth]{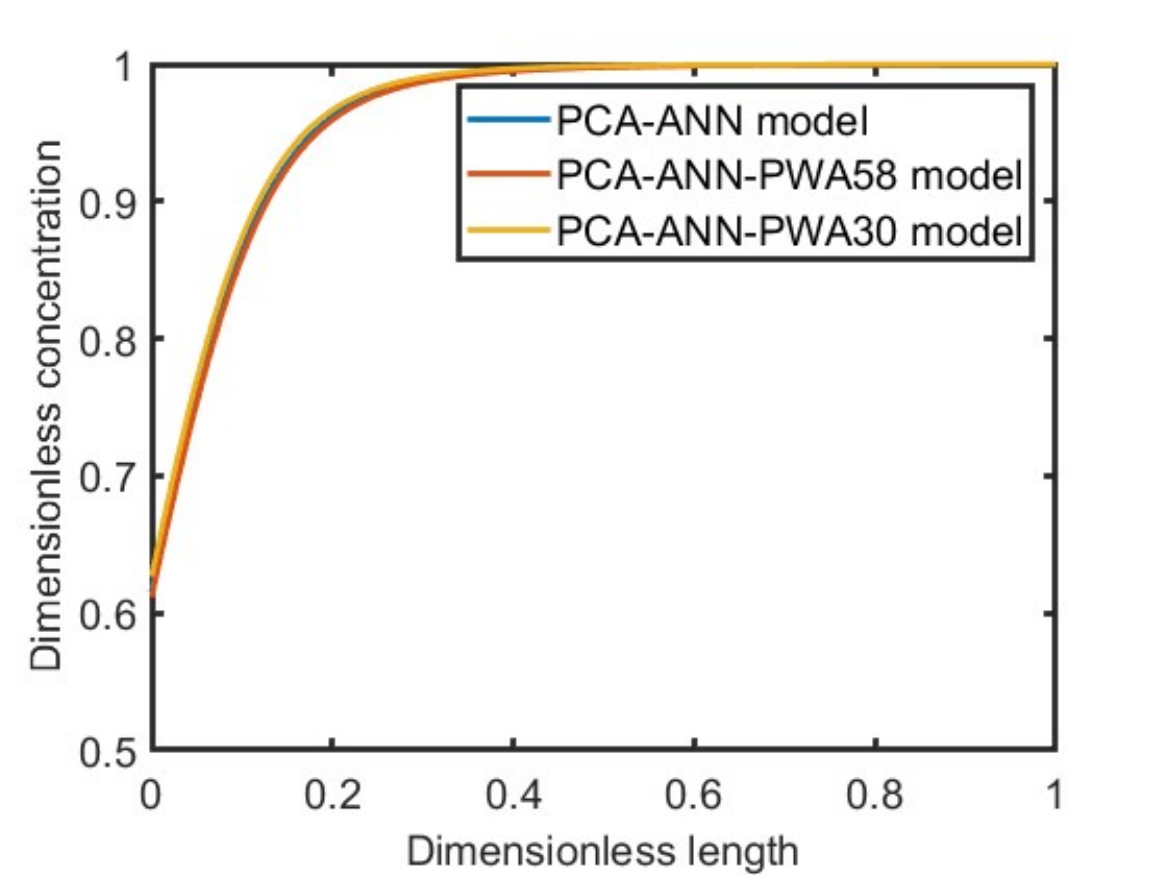}}
\caption{Solution profiles for dimensionless temperature and concentration (adaptive PWA)}
\label{solution_TR}
\end{figure}

\begin{figure}[H]
\centering
\includegraphics[width=0.5\textwidth]{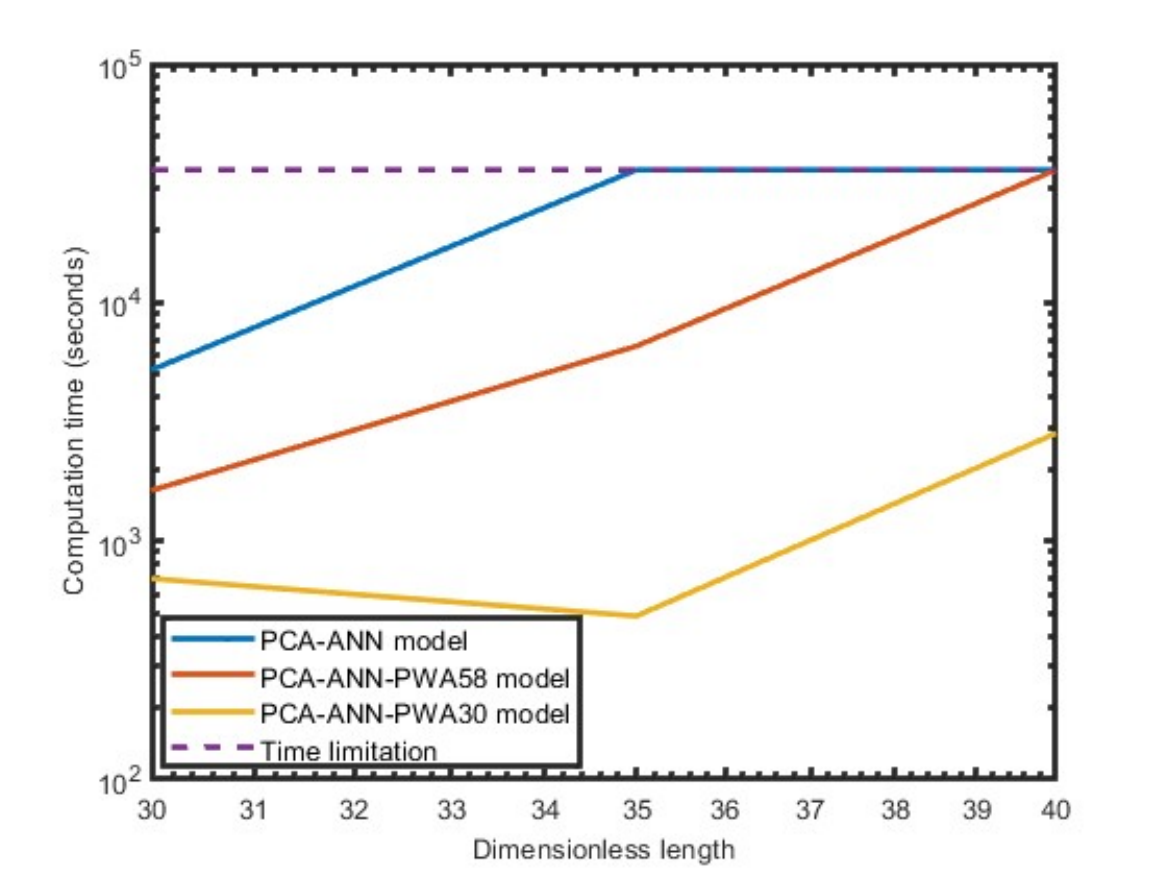}
\caption{Computational time (seconds) for different numbers of neurons (adaptive PWA)}
\label{solution_TR1}
\end{figure}

 Similarly, the above rector model was discretized in 250 central finite differences , which comprises our in-house full order model (FOM) simulator. PCA reduction was performed on 998 systematically collected LHC samples through the FOM first to reduce the 500 state variables down to 12. Subsequently 30- to 40-neuron ANNs were employed to obtain reduced PCA-ANN models comprising 3 inputs ( $T_{wi}$ ), and 12 state variables (outputs).  Two different PWA schemes (with 30 and 58 linear segments, respectively) were constructed following the above adaptive procedure. Here, PCA-ANN models with/without PWA steps formulates the reduced order models (ROMs).
 The optimisation results, depicted in Tb.\ref{Tb.2} and Fig.\ref{solution_TR},  are computed to compare the optimisation performance using the PCA-ANN and PCA-ANN-PWA models with 30 and 58 linear segments, respectively. In Tb.\ref{Tb.2}, "function value" refers to the function value through the FOM at the obtained optimal solutions of the surrogate models while "Global optimal value" represents the function value computed by performing direct global optimisation on the FOM  (solving the discretised algebraic equation formulation through global optimisation solver BARON).  All three computational cases converge to almost the same solutions, with the corresponding objective function value through the FOM close to the global optimal value of FOM .  The maximum error is 0.026\%, and the optimal solution profiles for concentration and temperature distributions (adaptive PWA scheme) are very close to each other for all models (Fig.\ref{solution_TR}). Fig.\ref{solution_TR1} compares the computational times required to perform optimisation using  the PCA- ANN and the PCA-ANN-PWA models with different number of ANN neurons. The limit time (max time for computations to stop) was set to be 36000 seconds. Both of the relative and absolute tolerances were set to be 0.0002.  The computational times increase rapidly with more neurons for all three kinds of surrogate models. Specifically, the computational costs reach the limit time for the PCA-ANN models with 35 and 40 neurons while the computations for the  PCA-ANN-PWA (30 linear segments) models completed with less than 7000 CPU seconds. It can be also seen that the computational time required is significantly less for PCA-ANN-PWA models, irrespective of the number of ANN neurons, implying the high computational efficiency of the proposed PWA methodology.  
\begin{table}[H]
\centering
  \caption{Optimal results (through full models) comparisons for surrogate tubular reactor models}
    \label{Tb.2}
  \smallskip
    \begin{tabular}{c|c|c|c|c|c}
  \hline
    Model &Solver&Function  &True&Relative&CPU\\
  (PCA-ANN)&& value&optimal value&error&(s)\\
  \hline
  (1 layer, 30 neurons, tanh)&BARON&0.99982&0.99998&0.016\%&5222.94\\
  PWA(30 adaptive linear segments)&CPLEX&0.99972&0.99998&0.026\%&695.63\\
  PWA(58 adaptive linear segments)&CPLEX&0.99986&0.99998&0.012\%&1631.14\\
  \hline
  \end{tabular}
\end{table}

Tab.\ref{Tb.2_1} shows the comparison of the optimal solutions computed using the reduced model (PCA-ANN model, a layer, 30 neurons, tanh) with adaptive and uniform linear segments. Almost the same optimal solutions are computed, which are close to the global optimum value of the FOM (0.99998). The proposed ANN-PWA model with 30 linear segments (adaptive linear scheme) required significantly less computational time, around one-fold reduction compared to the ANN-PWA model with 58 linear segments.  Furthermore, the uniform PWA formulation produced worse solutions than the proposed adaptive scheme, which could be explained by the large model predictive errors uniform PWA ( around 3 \% relative errors for all collected samples ) than that (0.005\%) of our adaptive scheme. Meanwhile, solving the same-size of the original MILP problem resulting from the uniform partitioning scheme is much faster than that of the proposed adaptive PWA scheme, less than 1s !  This  indicates the size of the MILP problems may be significantly reduced and/or much stronger relaxation tightness through automatically exploiting the special structure of reformulated problem using CPLEX solver. 
\begin{table}[H]
\centering
  \caption{Optimal result comparisons for adaptive and uniform PWA schemes}
    \label{Tb.2_1}
  \smallskip
    \begin{tabular}{c|c|c|c|c|c}
  \hline
    Model &Solver&Function & Global &Relative&CPU\\
  (PCA-ANN, 1 layer, 30 neurons, tanh)&&value&optimal value&error&(s)\\
  \hline
  PWA(30  adaptive linear segments)&CPLEX&0.99972&0.99998&0.026\%&695.63\\
  PWA(30  uniform linear segments)&CPLEX&0.99453&0.99998&0.545\%& $\leq$ 1 \\
  PWA(58  adaptive linear segments)&CPLEX&0.99986&0.99998&0.012\%&1631.14\\
  PWA(58  uniform linear segments)&CPLEX&0.99824&0.99998&0.174\%& $\leq$ 1\\
  \hline
  \end{tabular}
\end{table}
To verify the increased computational efficiency of the DNN models, global optimisation is first performed using four surrogate models (ANN, ANN-PWA, tanh-DNN and relu-DNN model, respectively) for the tubular reactor case presented above. Tab.\ref{Tb.2_2} shows the optimisation results. The small-scale tanh-DNN could replace the larger shallow ANN model, resulting in significant computational savings, of more than one order of magnitude. The relu-DNN model requires more neurons than the tanh-DNN model, due to the simpler structure of the relu activation function. Despite the fact that the relu-DNN model is larger, its optimisation cost is much lower than that of the tanh-NN models. The rapid global optimisation computations when using the relu-NN model are attributed to
the advanced MILP solver algorithm utilised. Furthermore, the computation cost using the relu-DNN model is much less than that that of using the ANN-PWA model with 30 linear segments because of the large(r) number of linear models involved in the PWA formulation. More linear models lead to more binary variables, requiring more branching steps hence reducing the computational efficiency. Meanwhile, the three PCA-DNN-PWA capture the wrong optimal solutions, where the relative errors can be around 80\%,  in despite of much fast computational speed for both adaptive and uniform schemes.  PWA errors can be propagated through the DNN structure and cause inaccurate PCA-DNN-PWA models, misleading the optimisation process although the structure of the "bad" models can be exploited to accelerate the computations. More linear segments are required to improve the accuracy of PCA-DNN-PWA models for further optimisation tasks. 
\begin{table}[H]
\centering
  \caption{Optimal result comparisons for shallow and deep PCA-ANN and/or PWA model}
    \label{Tb.2_2}
  \smallskip
    \begin{tabular}{c|c|c|c|c|c}
  \hline
    Model &Solver&Function&Global &Relative&CPU\\
  (PCA model)&&value&optimal value&error&(s)\\
  \hline
  ANN(1 layer, 30 neurons, tanh)&BARON&0.99982&0.99998&0.016\%&5222.94\\
  ANN-PWA(30  adaptive  segments)&CPLEX&0.99972&0.99998&0.026\%&695.63\\
  DNN(two layers, 8, 8 neurons, tanh)&BARON&0.99994&0.99998&0.004\%&109.95\\
  DNN-PWA(30, adaptive segments)&CPLEX&0.18376&0.99998&81.620\%&$\leq$ 1 \\
  DNN-PWA(58, adaptive segments)&CPLEX&0.25092&0.99998&74.910\%& $\leq$ 1\\
  DNN-PWA(58, uniform segments)&CPLEX&0.19537&0.99998&80.460\%& $\leq$ 1\\
  DNN(two layers, 40, 40 neurons, relu)&CPLEX&0.99976&0.99998&0.022\%&  42.22\\
  \hline
  \end{tabular}
\end{table}

To verify the increased computational efficiency of the PCA steps, the optimisation results of PCA-DNN double models and single DNN models are displayed as the below Tab.\ref{Tb.2_3}. It can be seen that the PCA-DNN models require less layers and neurons to account for the huge distributed information, resulting in  much less global optimisation cost, more than two orders of magnitude less than that of the direct DNN models for both tanh and relu activation functions.  Meanwhile, computations of ReLU based DNN models could converge to the optimal solutions much faster than that of tanh-based DNN since the advanced MILP solver algorithm is utilised.  

\begin{table}[H]
\centering
  \caption{Optimal result comparisons for DNN models with and without PCA steps}
    \label{Tb.2_3}
  \smallskip
    \begin{tabular}{c|c|c|c|c|c}
  \hline
    Model &Solver&Function & Global&Relative&CPU\\
  &&value&optimal value&error&(s)\\
  \hline
  PCA-DNN(two layers,&BARON&0.99994&0.99998&0.004\%&  109.95\\
   8, 8 neurons, tanh)&&&&& \\
   PCA-DNN(two layers, &CPLEX&0.99976&0.99998&0.022\%& 42.22\\
   40, 40 neurons, relu)&&&&& \\
   DNN(three layers,&BARON&-&0.99998&-\%&  $\geq$36000\\
   40, 40, 40 neurons, tanh)&&&&& \\
   DNN(three layers, &CPLEX&0.99994&0.99998&0.004\%& 13505.36\\
   50, 50, 50 neurons, relu)&&&&& \\
  \hline
  \end{tabular}
\end{table}

\section{Case study 2: combustion process}
To further test the significant advantages of PCA reduction and the deep rectifier neural network in our global optimisation formulation, a more challenging combustion process \cite{wei2012optimization,lang2009reduced} is considered here.
%
The combustion process takes place in a horizontal cylindrical combustor, 1.8m in length and 0.45m in diameter with a fuel nozzle that has a diameter 0.0045m. .The chemical reactions in the combustor are the following:
\begin{itemize}
    \item \ce{CH4 + 2O2 $\rightarrow$  CO2 + 2H2O}
    \item \ce{C2H4 + 3O2 $\rightarrow$  2CO2 + 2H2O}
    \item \ce{C3H8 + 5O2 $\rightarrow$  3CO2 + 4H2O}
    \item \ce{C4H10 + 8.5O2 $\rightarrow$  4CO2 + 5H2O}
\end{itemize}
In addition a complex \ce{NO} mechanism, comprising thermal \ce{NO}, prompt \ce{NO} and \ce{N2O} intermediate mechanisms is also taken into account. The fuel \ce{NO} mechanism was ignored due to the small amount of nitrogen in the feed. 
Thermal efficiency can be improved by increasing combustion temperature, which  however, inevitably leads to more pollutant emissions, such as \ce{NOx}. The \ce{NOx} production is dominated by the thermal \ce{NO} mechanism, given below, which is very sensitive to temperature.
\begin{itemize}
    \item \ce{O + N2 $\rightleftharpoons$  NO + N}
    \item \ce{O2 + N $\rightleftharpoons$  NO + O}
    \item \ce{N + OH $\rightleftharpoons$  NO + H}
\end{itemize}

This work focuses on the optimisation of inlet operational conditions (shown in Tb.\ref{Tb.4})
 
in order to minimize \ce {NOx} emissions. In addition to chemical reactions, multiple physical phenomena are involved, including complex turbulent flows, heat and mass transfer. Commercial CFD software was used, namely ANSYS/FLUENT, to construct high-fidelity CFD models to compute  velocity, temperature and component fraction fields.
\subsubsection{CFD Model Description}
The computational domain for the CFD model consisted of a 2-dimensional axisymmetric depiction of the combustor  To ensure that computations are grid independent, numerical experiments using 5481, 6381, 9081 and 14832 computational cells  were performed 
Finally, 9081 computational cells (9332 nodes) were chosen  as solutions did not change with more computational cells/nodes. The renormalisation group (RNG) $k-\epsilon$ turbulence model for fluid flow was employed. The eddy-dissipation model was employed for the species transport equations because the overall reaction rate was controlled by turbulent mixing. To take into account the  effects of thermal radiation, including absorption and scatting coefficients, a discrete ordinates (DO) radiation model was used. 

The second-order upwind scheme was applied for the space derivatives of the advection terms in all transport equations. The SIMPLE algorithm was employed to handle the velocity-pressure coupling in the flow field equations. Convergence criteria required the residual for the energy equation to be below $1\times 10^{-6}$  and residuals for the other model equations to be below $1\times 10^{-3}$. The mass-weighted-averages of temperature at the exit and the maximum temperature of the entire fluid were also monitored as other convergence criteria. 
 The base case inlet conditions used are given below. For the fuel gas, the base value for the inlet velocity was 100 m/s and that for the inlet temperature was 298K. The inlet composition was as follows: \ce{CH4}: 87.8\%, \ce{C2H4}: 4.6\%, \ce{C3H8}: 1.6\%, \ce{C4H10}: 0.5\%, \ce{N2}: 5.5\%. For preheated air, the inlet velocity was 85 m/s and the inlet temperature 1473 K while the  inlet composition was: \ce{O2}: 19.5\%, \ce{N2}: 59.1\%, \ce{H2O}: 15\%, \ce{CO2}: 6.4\%, \ce{NO}: 110 ppm.  Five independent variables were used to optimise the whole process. The independent variables along with their allowable ranges are listed in Tab.\ref{Tb.4}.
\begin{table}[H]
\centering
  \caption{Range of independent variables}
   \label{Tb.4}
  \smallskip
    \begin{tabular}{c|c|c}
  \hline
  Variables&Range&Units\\
  \hline
inlet air velocity&[85,125]&m/s\\
inlet fuel velocity &[80,120]&m/s\\
oxygen mass fraction (inlet air)& [18.5, 19.5] &\%\\
inlet air temperature & [1450, 1600] &K\\ 
inlet  fuel temperature& [298, 398] &K\\ 
  \hline
  \end{tabular}
\end{table}
\begin{figure}[H]
\centering
\includegraphics[width= 0.8\textwidth]{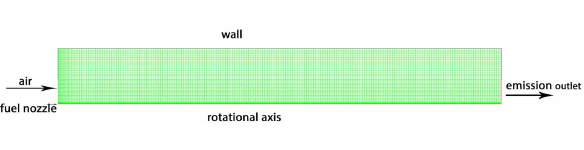}
\caption{Two-dimensional geometry of a single  axisymmetric combustor can and its mesh}
\label{fig:9}
\end{figure}
\subsubsection{Model reduction}
\par Although the high-fidelity CFD model can provide accurate simulation results,  its black-box characteristics and overall complexity make further optimisation and control tasks computationally tedious. Reduced surrogate models need to be developed to deal with the challenges arising. The LHC  sampling method was utilised to collect 1024 CFD samples, which took around 6 days using 4-CPU parallel computing. The input variables are the ones listed in Tab.\ref{Tb.4}, while the output results are the physical field data along with the average NOx emission at the outlet surface.  Due to the high dimensionality of the FOM, direct mapping of the input-output relationship would result in extremely large-scale ANN surrogate models, which could easily exceed the capability of current optimisation algorithms. Therefore, the PCA step was first employed and then surrogate ANN models were constructed based on the PCA-reduced models. 
%
ANN models were built for the field data, to construct the redcuced PCA-ANN constraints and for the average output NOx emission to formulate the ANN-reduced objective function. 
The field data 
include axial and radial velocity, temperature, \ce{N2}, \ce{H2O}, \ce{O2}, \ce{CO2}, \ce{C4H10}, \ce{CH4}, \ce{C2H4}, \ce{C3H8}, and \ce{NO} fraction concentrations (12 state variables). It should be noted that the average output NOx emission is only one variable so does not require PCA reduction. In this work, PCA was performed separately for each state variable.  While some PCA methods  compute principal components for all state variables together, we found that working on each state variable we could generate more accurate principal components. The standard criterion, of capturing 99.99\% of the total energy, was set. This way, the reduced surrogate models,  were built, as displayed in Tab.\ref{Tb.5}. 
\begin{table}[H]
  \caption{Number of PCs and corresponding ANN models}
  \label{Tb.5}
  \smallskip
  \centering
  \begin{tabular}{c|c|c|c}
  \hline
  Variables&No of PCs&DNN (2 layers, tanh)&DNN (2 layers, relu)\\
  &&No of neurons&No of neurons\\
  \hline
   Axial velocity&4& 14, 14&14, 14\\
Radial velocity&9& 15, 15&22, 22\\
  Temperature&6& 16, 16 &16, 16\\
 \ce{N2} concentration fraction&7& 19, 19&24, 24\\ 
 \ce{H2O} concentration fraction&8& 15, 15 &18, 18\\ 
 \ce{O2} concentration fraction&6& 17, 17 &20, 20\\
 \ce{CO2} concentration  fraction&7& 12, 12&14, 14\\ 
 \ce{C4H10} concentration fraction&6& 15, 15 &17, 17\\ 
 \ce{CH4} concentration fraction&7&26, 26 &28, 28\\ 
  \ce{C2H4}concentration fraction&6& 10, 10&18, 18\\ 
 \ce{C3H8}concentration fraction&6& 18, 18 &24, 24\\ 
 \ce{NO} concentration fraction&4&12, 12 &14, 14\\
   \hline
  Objective: output \ce{NOx} emission&-&ANN (1 layer, tanh)&ANN (1 layer, relu)\\
  &-&14 &30\\ 
  \hline
  \end{tabular}
\end{table}
\subsubsection{Model validation}
\begin{figure}[H]
\centering
\subfigure [Velocity (m/s) comparison among full model and reduced models]{
\label{fig:10}
\includegraphics[width= 0.45\textwidth]{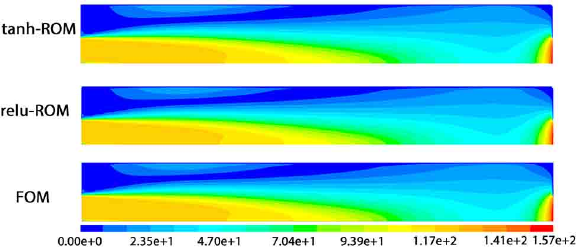}}
\hfill
\subfigure[Temperature (K) comparison among full model and reduced models]{
\label{fig:11}
\includegraphics[width= 0.45\textwidth]{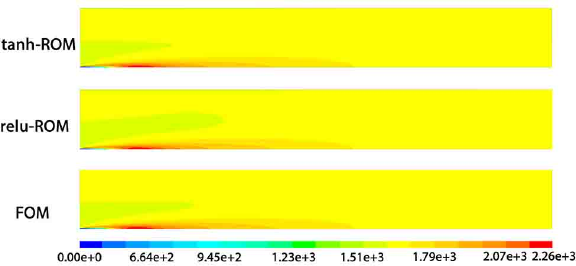}}
\subfigure[\ce{H2O} mass comparison among full model and reduced models]{
\label{fig:12}
\includegraphics[width= 0.45\textwidth]{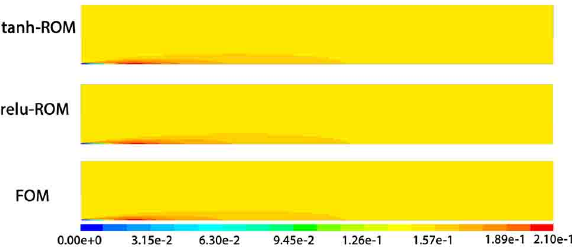}}
\hfill
\subfigure[\ce{O2} mass comparison among full model and reduced models]{
\label{fig:13}
\includegraphics[width= 0.45\textwidth]{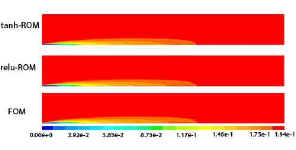}}
\subfigure[\ce{CO2} mass comparison among full model and reduced models]{
\label{fig:14}
\includegraphics[width= 0.45\textwidth]{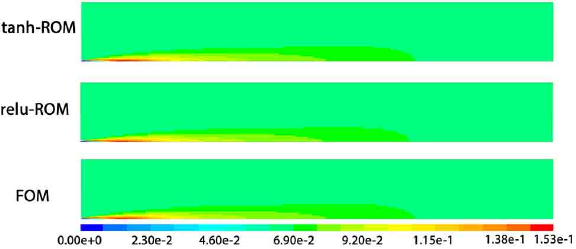}}
\caption{Comparison of velocity, temperature and concentration fraction field between FOM and ROMs}
\end{figure}

Model validation was performed for the reduced models before the subsequent optimisation step,  taking into account two aspects, representation ability and prediction ability. The representation ability of the reduced models was tested through the comparison between the FOM and the ROMs on the base case inlet conditions. Computational results show only very small  differences, especially for \ce{N2}, \ce{C4H10}, \ce{CH4}, \ce{C2H4}, \ce{C3H8}, and \ce{NO} fraction fields. The above species fraction fields are close to  uniform distribution across the combustor, except for the small area near the fuel nozzle. Fig.\ref{fig:10}, \ref{fig:11}, \ref{fig:12}, \ref{fig:13}, \ref{fig:14} depict the velocity field, temperature field, \ce{O2} , \ce{CO2} and \ce{H2O} concentration fraction field  of FOM, tanh-ROM and relu-ROM, respectively for the inlet base values. The five contour diagrams illustrate that flow, temperature and mass fraction fields of FOM and ROMs are very close, indicating the strong representation ability of the ROMs. Moreover, the tanh-DNN reduced models show smaller difference from the FOM than the relu-DNN reduced models, especially for the temperature field, implying the better accuracy of the tanh-DNN models due to the non-linearity of the tanh function. 
Tab.\ref{Tb.6} shows the comparison of maximum field values between FOM and ROMs and the corresponding errors. The largest error is only 0.56\%. 
To test the ROMs prediction ability, 40 random inlet condition points different than the base case ones were chosen and compered with FOM results. The largest error was less than 5\%   indicating that the ROMs can be reliably used for further optimisation studies.
Fuerhtermore, the ROMs exhibit significant computational savings compared to the full-order CFD models as expected. 
The average CPU time for the CFD model (run in ANSYS/FLUENT) is approximately 1560 CPU seconds, while each ROM  requires less than 0.1 CPU seconds and can be efficiently used to perform global optimisation studies. 
\begin{table}[H]
   \caption{Average value comparison of FOM and ROMs }
   \label{Tb.6}
  \smallskip
  \centering
  \begin{tabular}{c|c|c|c|c|c}
  \hline
  Variables&FOM&relu-ROMs&errors&tanh-ROMs&errors\\
  \hline
   Velocity(m/s)&29.82089& 29.6652&0.40\%&29.70166&0.56\%\\
  Temperature(K)&1625.259& 1621.948&0.03\%& 1625.702&0.20\%\\
 \ce{H2O} mass fraction&0.151743& 0.1518191 &0.02\%&0.1517722&0.05\%\\
 \ce{O2} mass fraction&0.1906906& 0.1908453&0.01\%&0.1907176&0.08\%\\
 \ce{CO2} mass fraction&0.0663304& 0.06627689&0.00\%&0.06633245&0.08\%\\
 Computational time for each sample&1560& $\leq$0.1&-&$\leq$0.1&-\\
  \hline
  \end{tabular}
\end{table}
\subsubsection{Global optimisation}
\par In this section, global optimisation is implemented using the validated reduced models. The general mathematical formulation is given in
Eq.\eqref{eq25}. In the combustion optimisation problem, $\boldsymbol{d}$ are the 5 inlet operation parameters, and $\boldsymbol {u'}$ are the 76 reduced state variables . The objective function $G(\boldsymbol {d},\boldsymbol{u'})$ represents the average outlet \ce{NOx} emission. 
The allowable ranges for the input variables are given in Tab.\ref{Tb.4}, while the bounds for the state variables are given in Tab.\ref{Tb.7}. 
It should be noted that the state variable bounds are implemented through the inverse projection
\begin{equation}\label{inv1}
\begin{aligned}
lb\leq\boldsymbol{P^T}\boldsymbol{u'} \leq ub \\
  \end{aligned}
\end{equation}
where $lb$ and $ub$ denote lower and upper bounds, respectively. 

Finally, a MILP problem with 29,903 linear constraints, corresponding to the equality constraints and 488 binary variables corresponding to the total number of ANN neurons is formulated for the relu-based ROM, while a NLP problem with 28247 linear constraints, 
and 392 nonlinear terms is constructed for the tanh-based ROM.  
The limit value for the computational time was set to be 100 hours. Both of the relative and absolute tolerances were set to be 0.002.  

The NLP problem did not converge to a feasible solution in BARON within the allowable time, probably due to the high non-convex activation function $tanh$ and large number of variables than inhibited the branch-and-bound algorithm. 

\begin{figure}[H]
\centering
\subfigure[Velocity (m/s) comparison among full model and reduced models under optimal condition]{
\label{fig:15}
\includegraphics[width= 0.45\textwidth]{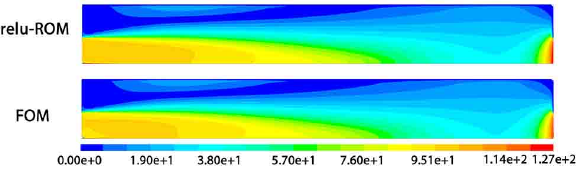}}
\subfigure[Temperature (K) comparison among full model and reduced models under optimal condition]{
\label{fig:16}
\includegraphics[width= 0.45\textwidth]{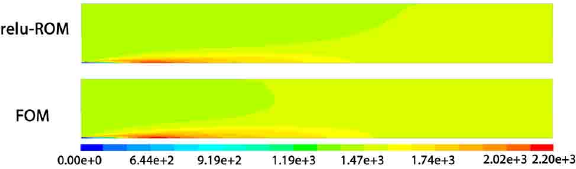}}
\hfill
\subfigure[\ce{H2O} mass fraction comparison among full model and reduced models under optimal condition]{
\label{fig:17}
\includegraphics[width= 0.45\textwidth]{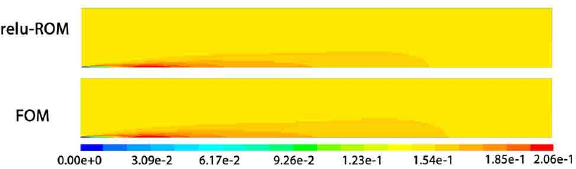}}
\subfigure[\ce{O2} mass fraction comparison among full model and reduced models under optimal condition]{
\label{fig:18}
\includegraphics[width= 0.45\textwidth]{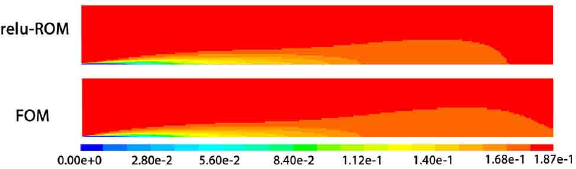}}
\hfill
\subfigure[\ce{CO2} mass fraction comparison among full model and reduced models under optimal condition]{
\label{fig:19}
\includegraphics[width= 0.45\textwidth]{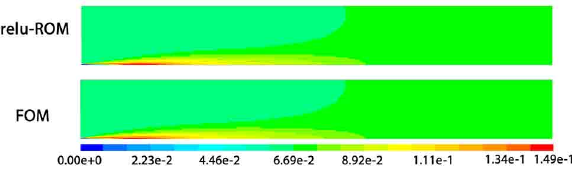}}
\caption{Comparison of optimal velocity, temperature and fraction concentration field between FOM and ROM}
\end{figure}

The $relu$-based MILP problem converged in 501.89s in CPLEX. The computed optimal solution was: NOx emission: 110.17 ppm, air velocity: 95.07 m/s, fuel velocity:119.08m/s, oxygen fraction concentration (air): 18.50 \%,  air temperature: 1450 K and  fuel temperature: 369.83K. 
To validate the computed optimal solutions, we performed a full CFD simulation in ANSYS/FLUENt using the calculated optimal inlet conditions. 
The outlet \ce{NOx} emission was computed to be 113.26 ppm, which was very close to the calculated optimum with an error of approximately 2.73 \%, which is small enough for most industrial cases. Fig.\ref{fig:15}, \ref{fig:16}, \ref{fig:17}, \ref{fig:18} ,\ref{fig:19}, depict a comparison of the main field state varibales at the optimal; conditions computed by the reduced and the full models, respectively. As it can be observed, the optimal solution computed through the ROM is very close to FOM simulation using the optimal inlet conditions. Tab.\ref{Tb.8} gives a comparison of the corresponding max values across the whole domain. The  performance of the reduced model is very close to the full model with the biggest error being less than 3\%. The computational cost for the ReLu-based MILP problem is significantly reduced compared to the NLP problem, which in this case could not converge, signifying the efficiency of the our model reduction-based global optimisation methodology.   
\begin{table}[H]
  \caption{Range of state variables }
   \label{Tb.7}
     \smallskip
     \centering
     \begin{tabular}{c|c|c}
  \hline
 State variables&Range&Units\\
  \hline
   Axial velocity&[-150,150]&m/s\\
 Radial velocity &[-150,150]&m/s\\
  Temperature & [0, 2200] &K\\
\ce{N2} concentration fraction& [0,1] &-\\ 
 \ce{H2O} concentration fraction& [0,1] &-\\ 
  \ce{O2} concentration fraction& [0,1] &-\\ 
   \ce{CO2} concentration fraction& [0,1] &-\\ 
    \ce{C4H10} concentration fraction& [0,1] &-\\ 
    \ce{CH4} concentration fraction& [0,1] &-\\ 
    \ce{C2H4} concentration fraction& [0,1] &-\\ 
    \ce{C3H8} concentration fraction& [0,1] &-\\ 
    \ce{NO} concentration fraction& [0,1] &-\\ 
  \hline
  \end{tabular}
  \end{table}  

 \begin{table}[H]
   \caption{Average value comparison of FOM and ROMs }
  \label{Tb.8}
  \smallskip
  \centering
  \begin{tabular}{c|c|c|c|c|c}
  \hline
 Variables&FOM&relu-ROMs&errors\\
  \hline
   Velocity(m/s)&23.11901& 23.73662&2.67\%\\
  Temperature(K)&1487.854& 1456.745&2.09\%\\
 \ce{H2O} mass fraction&0.1522371& 0.1520701&0.11\%\\
 \ce{O2} mass fraction&0.1803872& 0.1806236&0.13\%\\
 \ce{CO2} mass fraction&0.06699268& 0.06689615&0.14\%\\
  Output \ce{NOx} emission (ppm)&113.26& 110.17&2.73\%\\
  \hline
  \end{tabular}
\end{table}
 \section{Conclusions}
\par This paper presents  model-reduction, machine-learning based global optimisation framework for large-scale nonlinear steady-state systems. A double model reduction, comprising principal component analysis and artificial neural networks, were first employed  to construct accurate surrogate reduced models through a machine-learning step, which was then utilised by deterministic global optimisation methods. The high non-convexity of  the activation function in reduced ANN models affects the computational speed branch-and-bound algorithms. To overcome this barrier, two improvements were proposed.  Firstly, when nonlinear activation functions are preferable, a piece-wise affine reformulation to transform the nonlinear branching into binary variables resulting in a MILP problem with higher computational efficiency. Secondly, the implementation of a continuous piece-wise linear activation function-based deep ANN structure to improve computational accuracy. A number of case studies of different size including  a tubular reactor and a complex large-scale combustion process were employed to illustrate the favorable performance of the presented framework. 
Despite the efficient implementation of the developed methodology, there are still remaining challenges to efficiently compute the global optimum for large-scale optimisation problems. 
Firstly, this work assumes enough representative samples as a basis to construct the reduced order models. Smart sampling methods to achieve optimal trade-off between quality and quantity are important for improving both efficiency and accuracy, as well as verification methods to guarantee the accuracy of the computed solutions\cite{botoeva2020efficient}.
Secondly, global optimisation even using reduced surrogate models is still computationally expensive. Advanced data techniques and MILP algorithms \cite{anderson2018strong} may further improve computational efficiency of this optimisation framework.

\section*{acknowledgements}
The financial support of the University of Manchester and China Scholarship Council joint scholarship for MT's PhD studies is gratefully acknowledged 
\section*{conflict of interest}
The authors declare there is no conflict of interest. 

\bibliographystyle{model1-num-names}
\bibliography{sample.bib}
\end{document}